\documentclass[a4paper,12pt]{article}
\usepackage[warn]{mathtext}
\usepackage{cmap}
\usepackage[T2A]{fontenc}
\usepackage[cp1251]{inputenc}
\usepackage[russian]{babel}
\usepackage{amsmath}
\usepackage{wasysym}
\usepackage{amssymb,amsthm}
\usepackage[dvips]{graphicx}

\theoremstyle{plain}
\newtheorem{theorem}{Теорема}[section]

\newtheorem{lemma}[theorem]{Лемма}
\newtheorem{proposition}[theorem]{Предложение}

\newtheorem{algoritm}[theorem]{Алгоритм}

\theoremstyle{definition}

\newtheorem{example}[theorem]{Пример}

\begin{document}

\title{Равномерная рекуррентность морфических слов}

\author{И.~Митрофанов \footnote{ работа выполнена при частичной финансовой поддержке фонда Дмитрия Зимина ``Династия" и гранта РФФИ № 14-01-00548}
}
\date{}
\maketitle

\section{Введение.}
Бесконечное вправо слово $W$ над конечным алфавитом называется {\it равномерно рекуррентным}, если для любого конечного подслова $u$ существует такое $k\in \mathbb{N}$,
что $u$ встречается в любом подслове слова $W$ длины $k$.

Понятие равномерной рекуррентности возникло из динамических систем (равномерная рекуррентность слова означает, что порождённый им субшифт минимален), 
также его можно рассматривать как одно из обобщений свойства периодичности бесконечного слова \cite{PS}.

В работах \cite{Pr} и \cite{PS} изучается равномерная рекуррентность (почти периодичность) {\it морфических слов} и ставится алгоритмическая задача проверки равномерной рекуррентности морфических слов:

{\bf Вход:} два конечных алфавита $A$ и $B$, подстановка $\varphi:A^*\to A^*$, продолжающаяся над буквой $a_1\in A$, кодирование $\psi:A^*\to B^*$.\\
{\bf Вопрос:} является ли слово $\psi(\varphi^{\infty}(a_1))$ равномерно рекуррентным?

В работе \cite{Pr} был найден алгоритм проверки равномерной рекуррентности для двух классов морфических слов: 
для {\it чисто морфических} (то есть для слов вида $\varphi^{\infty}(a_1)$)
и для {\it автоматных последовательностей}, то есть для подстановочных систем с условием 
$|\varphi(a_i)|=|\varphi(a_j)|$ для любых $i$ и $j$.

Также в \cite{Pr} и \cite{PS} была высказана гипотеза, что задача проверки равномерной рекуррентности алгоритмически разрешима для любого морфического слова.

Автор и F.Durand независимо получили доказательство этой гипотезы:

\begin{theorem}\label{UR} \cite{Mur,Dur}
Задача проверки равномерной рекуррентности морфических слов алгоритмически разрешима.
\end{theorem}

Доказательство автора опиралось на результат совместной с А.Я.Беловым статьи \cite{MKR}, позже было получено другое доказательство.

\section{План построения алгоритма.}

Сформулируем следующую алгоритмическую задачу, известную как HD0L-periodicity problem:

{\bf Вход:} два конечных алфавита $A$ и $B$, подстановка $\varphi:A^*\to A^*$, продолжающаяся над буквой $a_1\in A$, кодирование $\psi:A^*\to B^*$.\\
{\bf Вопрос:} является ли слово $\psi(\varphi^{\infty}(a_1))$ периодичным с некоторого момента?

Автор и F.Durand независимо установили, что
\begin{theorem}\cite{Dur2,M,MPF}\label{HD0L}
HD0L-periodicity problem алгоритмически разрешима.
\end{theorem}

В построении алгоритма важную роль играет свойство минимальностя языков равномерно-рекуррентных слов.

\begin{proposition}\label{minimal}
Пусть $W_1$ и $W_2$ -- два сверхслова такие, что $W_2$ равномерно рекуррентно 
и все конечные подслова $W_2$ являются подсловами $W$.
Тогда $W_1$ равномерно рекуррентно если и только если все его подслова являются подсловами $W_2$.
\end{proposition}

\subsection*{Замены подстановочной системы.}

Пусть нужно определить равномерную рекурренность сверхслова 
$W=\psi(\varphi^{\infty}(a_1))$.

Прежде всего, можно считать, что морфизм $\varphi$ нестирающий, т.е. $|\varphi(a)|>0$
для любой буквы $a\in A$.

Это следует из следующей теоремы:
\begin{theorem}[\cite{AS}, часть $7$] Пусть $f: A^*\to B^*$ и $g:A^*\to A^*$ --- произвольные морфизмы, и пусть $f(g^{\infty}(a_1))$ --- бесконечное вправо слово. 
Тогда существует алфавит $A'$, буква $a'_1\in A'$, нестирающая подстановка
$\varphi:A'^*\to A'^*$, продолжаемая над $a'_1$,
и кодирование $\tau:A'\to B$ такие, что $f(g^\infty(a_1))=\tau(\varphi^{\infty}(a'_1))$.
\end{theorem}

Более того, новые морфизмы можно найти алгоритмически.

Слово $w\in A^*$ называется $\varphi-${\it ограниченным}, если последовательность длин
$$
w,\varphi(w),\varphi^2(w),\varphi^3(w),\dots
$$
с какого-то места зацикливается.

Иначе $|\varphi^n(w)|\rightarrow \infty$ при $n\rightarrow \infty$ и $w$ называется {\it $\varphi-$растущим}.
Конечное слово является $\varphi-$ограниченным если и только если оно состоит из $\varphi-$ограниченных символов.

\begin{lemma} \label{finwords}
Следующая задача алгоритмически разрешима:

{\bf Вход:} подстановка $\sigma $ и буква $a$.

{\bf Выход:} 
1. Ответ на вопрос <<верно ли, что набор $\sigma -$ограниченных подслов сверхслова $\sigma ^{\infty}(a)$ бесконечен?>>\\
2. Если ответ на первый вопрос <<да>>, то такое конечное непустое слово $U$, что 
$\sigma ^{\infty}(a)$ содержит любую степень этого слова $U^k$.\\
3. Если ответ на первый вопрос <<нет>>, то список всех $\sigma -$ограниченных подслов $\sigma ^{\infty}(a)$.
\end{lemma}

Применим алгоритм из этой леммы к $\varphi$ и $a_1$. 
Дальнейшее зависит от того, какой ответ получен на первый вопрос.

{\bf 1. Набор бесконечен.} 
Тогда любое конечное подслово слова периодичного сверхслова $(\psi(U))^{\infty}$ является подсловом $W$.
Так как чисто периодичное слово является равномерно рекуррентным, то, 
согласно лемме \ref{minimal}, $W$ может быть равномерно рекуррентно в том и только том случае, если оно периодично.
По теореме \ref{HD0L} периодичность $W$ алгоритмически определяется, следовательно, определяется и равномерная рекуррентность. 

{\bf 2. Набор конечен.} В этом случае можно найти ещё более <<хорошую>> подстановочную систему, порождающую такое же сверхслово.

\begin{lemma}\cite{Pr,Pans}
Пусть в $\varphi^{\infty}(a_1)$ конечное число $\varphi$-ограниченных подслов.
Тогда можно построить конечный алфавит $C$, подстановку $g:C^*\to C^*$, 
продолжающуюся над буквой $c_1\in C$ и нестирающий морфизм $h:C^*\to A$ такие, что\\
1. $\varphi^{\infty}(a_1)=h(g^{\infty}(c_1))$;\\
2. Все буквы алфавита $C$ являются $g$-растущими.
\end{lemma}

Таким образом, задача \ref{UR} сводится к такой задаче:\\
{\bf Вход:} подстановочная система $(A,a_1,\varphi,B,\psi)$, где 
все буквы алфавита $A$ являются $\varphi$-растущими,
а $\psi$ -- нестирающий морфизм.\\
{\bf Вопрос:} является ли слово $\psi(\varphi^{\infty}(a_1))$ равномерно рекуррентным?

Назовём подстановку $\sigma$, действующую на алфавите $S$, {\it примитивной}, если существует такое натуральное $k$, что для любых двух букв $s_i,s_j\in S$ выполнено
$s_i\in \sigma^{k}(s_j)$.

Широко известен следующий классический результат:
\begin{proposition} \cite{AS}
Если $\sigma:A^*\to A^*$ -- примитивная подстановка, 
а $h:A^*\to B^*$ -- произвольный нестирающий морфизм, то слово $h(\varphi^{\infty}(a_1))$
является равномерно рекуррентным. 
\end{proposition}

Далее мы будем сравнивать множество подслов $\psi(\varphi(a_1))$ и некоторого сверхслова, порождённого примитивной подстановкой.

\begin{lemma}\cite{Pr} \label{pr_ins}
Пусть все буквы алфавита $A$ являются $\varphi$- растущими. 
Тогда можно найти в алфавите $A$ подмножество $A'\subseteq A$ такое, 
что некоторая степень $\varphi$ является в ограничении на $A'$ примитивной подстановкой,
продолжающейся над $a'\in A'$.
\end{lemma}

Сверхслово $W'=\psi(\varphi^{\infty}(a'))$ является равномерно рекуррентным и все его слова являются подсловами сверхслова $W=\psi(\varphi^{\infty}(a_1))$.
Согласно \ref{minimal}, $W$ равномерно рекуррентно если и только если любое его конченое подслово является подсловом $W'$.

Согласно \ref{HD0L}, можно проверить периодичность слов $W$ и $W'$. 
Если $W'$ периодично, то равномерная рекуррентность $W$ равносильна периодичности $W$.
Далее полагаем, что $W'$ оказалось непериодично.

\subsection*{Исследование скоростей роста букв.}

\begin{proposition} Пусть $\sigma:A^*\to A^*$ --- произвольная подстановка. 
Тогда для любой буквы $a\in A$ выполняется одно из условий:\\
 $\exists k\in \mathbb N: \sigma^k(a)$ --- пустое слово.\\
Существуют числа 
$d(a)\in \mathbb N_0, c(a)\in \mathbb R_+, C(a)\in \mathbb R, \theta(a) \in \mathbb R$
такие, что  для всех натуральных $n$ выполнено
$$c(a)<\frac{\sigma^n(a)}{c(a)n^{d(a)}\theta(a)^n}<C(a).$$

\end{proposition}

Если подстановка $\varphi$ -- нестирающая, то для любой буквы $a\in A$ 
выполняется второе условие.
Пара $(d(a),\theta(a))$ называется {\it скоростью роста} буквы $a$. 

Скорости роста можно сравнивать:
 $(d_1,\theta_1)<(d_2,\theta_2)$ если $\theta_1<\theta_2$ или 
 $\theta_1=\theta_2$ и $d_1<d_2$.

Положим $(D,\Theta)=\max_{a_i\in A}(d(a_i),\theta(a_i))$. 
Все буквы $a$ со скоростью роста $(D,\Theta)$ назовём {\it $\varphi$-быстрорастущими.}
Понятно, что одной их $\varphi$-быстрорастущих букв является $a_1$.

\begin{lemma}
В сверхслове $W = \varphi^{\infty}(a_1)$ бесконечно много быстрорастущих букв.
\end{lemma}

\begin{lemma}
Если $W=\varphi^{\infty}(a_1)$ равномерно рекуррентно, то в $W$ 
не может встретиться подслово вида
$a_iUa_j$, где $a_i$ и $a_j$ -- $\varphi$-быстрорастущие буквы, а 
$U$ -- непустое слово без $\varphi$-быстрорастущих букв.
\end{lemma}

Отсюда следует, что либо все буквы $A$ имеют одинаковую скорость роста (при этом $D=0$), 
либо слово $W$ не является равномерно рекуррентным.

\begin{proposition}
Следующая задача алгоритмически разрешима:\\
{\bf Вход:} Подстановочная система $(A,a_1,\varphi,B,\psi)$\\
{\bf Вопрос:} Верно ли, что у всех букв $a_i\in A$ одинаковая скорость роста?
Если да, то найти такие числа $K_1$, $K_2$, $1<\Theta_1<\Theta_2$, что
для некоторого $\theta\in [\Theta_1;\Theta_2]$ и всех $a_i\in A$, $k\in \mathbb{N}$ выполнено $K_1\lambda^k<|\psi(\varphi^k(a_i))|<K_2\lambda^k$.
\end{proposition}

Таким образом, из лемм \ref{minimal} и \ref{pr_ins} следует, что задача о равномерной рекуррентности сводится к следующей алгоритмической задаче:

{\bf Дано:}\\
1. Три конечных алфавита $A$, $B$, $C$;\\
2. Числа $K_1$, $K_2$, $1<\Theta_1<\Theta_2$;\\
3. Четыре морфизма 
$\varphi:A^*\to A^*$, $\psi :A^*\to C^*$, $g:B^*\to B^*$, $h:B^*\to C^*$ такие, что
\begin{enumerate}
\item[a)] все морфизмы нестирающие;
\item[b)] морфизм $g$ продолжается над $b_1$, морфизм $\varphi$ примитивен и продолжается над $a_1\in A$;
\item[c)] для некоторого $\lambda\in [\Theta_1;\Theta_2]$ и всех $a_i\in A$, $b_j\in B$, $k\in \mathbb{N}$ выполнено $K_1\lambda^k<|\psi(\varphi^k(a_i))|<K_2\lambda^k$ и $K_1\lambda^k<|h(g^k(b_j))|<K_2\lambda^k$.
\item[d)] сверхслово $\psi(\varphi^{\infty}(a_1))$ непериодично;
\end{enumerate}
{\bf Определить:} верно ли, что любое подслово сверхслова $h(g^{\infty}(b_1))$ является подсловом сверхслова $\psi(\varphi^{\infty}(a_1))$?

\begin{theorem}\label{primitive}
Эта задача алгоритмически разрешима.
\end{theorem}

Дальнейшие конструкции строятся для доказательства теоремы \ref{primitive}.

\subsection*{Схемы расположения подслов.}

{\it Узлом} конечного слова $u$ длины $n$ будем называть одну из $n+1$ позиций: начало слова ({\it начальный} узел), конец слова ({\it конечный} узел) или один из $n-1$ промежутков между его буквами ({\it обычный} узел).

Каждые пара узлов определяет некоторое подслово в $u$, возможно, пустое.

Пусть $U=(U_1,U_2,\dots,U_n)$ и $u=(u_1,u_2,\dots,u_m)$ -- 
два упорядоченных множества конечных слов над одним и тем же алфавитом. 
Для каждого $U_i\in U$ определим множество {\it интересных} узлов: начальный, конечный, а также все те узлы, которые являются концами или началами каких-либо вхождений слов из $u$.

Назовём {\it схемой расположения подслов} для $U$ и $u$ следующую пару: 
\begin{enumerate}
\item Упорядоченное множество из $n$ чисел $t_1, t_2, \dots, t_n$, где $t_i$ -- это число интересных узлов в слове $U_i$.
\item Таблицу размером $n\times m$, в клетках которой находятся множества упорядоченных пар чисел, полученные по следующему правилу:
в клетке, стоящей на строке $i$ в столбце $j$ парами чисел, не превосходящих $t_i$, описываются все вхождения слова $u_j$ в $U_i$ (каждая пара чисел задаёт начало и конец некоторого вхождения).
\end{enumerate}

Такую схему будем обозначать $S(U,u)$.
Далее схемы вхождений подслов часто называются просто {\it схемами}.

\begin{example}\label{table_1}{} $n=3$, $m = 2$, $U_1 = abcabc$, $U_2 = bcabca$, $U_3 = bbbbbc$,
$u_1 = abc$, $u_2 = bc$.

В слове $U_1$ пять интересных узлов: ${}_1a_2bc_3a_4bc_5$. 

В слове $U_2$ также пять интересных узлов: ${}_1bc_2a_3bc_4a_5$.

А в слове $U_3$ три интересных узла: ${}_1bbbb_2bc_3$.

Тогда $S((U_1, U_2, U_3),(u_1, u_2))$ состоит из
вектора $(5,5,3)$ и следующей таблицы:

\begin{tabular}{|c|c|}
\hline 
(1,3), (3,5) & (2,3), (4,5) \\ 
\hline 
(2,4) & (1,2), (3,4) \\ 
\hline 
$\emptyset$ & (2,3) \\ 
\hline 
\end{tabular} 

\end{example}

Упорядоченное множество подслов слова $\varphi^{\infty}(a_1)$,
состоящих из одной или двух букв, назовём {\it порождающими словами} и будем обозначать 
$A_G$.
Все слова из $A_G$ алгоритмически находятся,
порядок на $G$ выбирается произвольный.

Введём обозначения: $U^k:=\{\psi(\varphi^k(x))|x\in A_G\}$,
$V^k:=\{h(g^k(x))|x\in B\}$.
Мы считаем, что элементы множества $U^k$ упорядочены так же, 
как соответствующие элементы $A_G$, а элементы $V^k$ упорядочены так же, как буквы в алфавите $B$.

\begin{proposition}\cite{AS}  \label{recurrence}
Существует такое $K(I)$, что если $u_1$ и $u_2$ --- два подслова сверхслова $\psi(\varphi^{\infty}(a_1))$
 такие, что $|u_2|>K(I)|u_1|$, то $u_1$ является подсловом $u_2$.
\end{proposition}

Запись $d(I)$ означает, что $d$ -- это число, 
алгоритмически определяемое по входным данным.

\begin{lemma}\label{S_per}
Можно алгоритмически найти число $d(I)$ такое, что из условий 
$k_1>l_1+d(I)$, $k_2>l_2+d(I)$ и $S(U^{k_1},V^{l_1})=S(U^{k_2},V^{l_2})$ следует 
$S(U^{k_1+1},V^{l_1+1})=S(U^{k_2+1},V^{l_2+1})$. 
Равенство схем вхождения подслов следует понимать как совпадение таблиц.
\end{lemma}

Доказательство основано на двух леммах:
\begin{lemma}\label{S_per_1}
$S(U^k,V^{l+1})$ можно однозначно определить по $S(U^k,V^l)$.
\end{lemma}

\begin{lemma}\label{S_per_2}
Пусть известна схема $S(U^k,V^l)$ и кроме этого известно, что 
$\min_i(|U^k|)>2 \max_j(|V_j|)$. 
Тогда можно алгоритмически найти $S(U^{k+1},V^l)$.
\end{lemma}

\subsection*{Алгоритм для задачи \ref{primitive}}.

{\it Размером} схемы вхождения подслов будем называть наибольшее количество элементов 
в клетках соответствующей таблицы. 
Так, размер схемы из примера \ref{table_1} равен 2.

Выберем $D(I)$ такое, что $D(I)>d(I)$ и $D(I)>\log_{\Theta_1}(\frac{K(I)K_2}{K_1})$.

Рассмотрим последовательность схем $S(U^{k+D(I)},V^k)$ при $k=1, 2, 3,\dots$

Положим $N(I)=\lceil \frac{2K_2}{K_1}K(I)\Theta_2^{D(I)}\rceil$.

Из непериодичности $\psi(\varphi^{\infty(a_1)})$ следует
\begin{lemma}\label{S_b}
При любом $k$ размер схемы $S(U^{k+D(I)},V^k)$ не превосходит $N(I)$.
\end{lemma}

Теперь можно привести псевдокод алгоритма, решающего задачу \ref{primitive}:

\begin{algoritm}\label{main_p}
\begin{enumerate}
\item Вычислить $D(I)$ и $N(I)$;
\item Последовательно для всех натуральных $k$ строить $S(U^{k+D(I)},V^{k})$, 
пока две из построенных схем не совпадут.
\item Если в двух совпавших схемах есть пустые клетки, то ответ <<нет>>, иначе ответ <<да>>.
\end{enumerate}
\end{algoritm}

Правильность работы алгоритма следует из лемм \ref{S_per} и \ref{recurrence}.

\section*{Благодарности}
Автор благодарит А.Я.Белова, А.В.Михалёва и Ф.Дюранда за внимание к работе.

\end{document}